\documentclass[dvips]{amsart}
\usepackage{math,epsfig,psfrag,bbold}

\newcommand\PSL{{\mathsf{PSL}}}

\def\pair<#1>{{\ll}#1{\gg}}
\def\9{{\mathbb 1}}
\def\8{{\mathbb 2}}

\begin{document}
\title{A Wilson Group of Non-Uniformly Exponential Growth}
\dedicatory{To Rostislav I. Grigorchuk on the occasion of his 50th birthday}
\date{\today}
\author{Laurent Bartholdi}
\address{Department of Mathematics, Evans Hall, U.C.Berkeley, USA}
\email{laurent@math.berkeley.edu}
\begin{abstract}
  This note constructs a finitely generated group $W$ whose
  word-growth is exponential, but for which the infimum of the growth
  rates over all finite generating sets is $1$ --- in other words, of
  non-uniformly exponential growth.
  
  This answers a question by Mikhael Gromov~\cite{gromov:metriques}.
  
  The construction also yields a group of intermediate growth $V$ that
  locally resembles $W$ in that (by changing the generating set of
  $W$) there are isomorphic balls of arbitrarily large radius in $V$
  and $W$'s Cayley graphs.
\end{abstract}
\maketitle

%%%%%%%%%%%%%%%%%%%%%%%%%%%%%%%%%%%%%%%%%%%%%%%%%%%%%%%%%%%%%%%%
\section{Introduction}
The purpose of this note is to construct in an as short and elementary
way as possible a group of non-uniformly exponential growth, i.e.\ a
group of exponential growth with a family of generating sets for which
the growth rate tends to $1$. The ``limit'' of these generating sets
generates a group of intermediate growth.

This construction is an adaptation of~\cite{bartholdi:interm}, which
describes a family of groups of intermediate growth. It came after
John Wilson announced he had produced such a group; my method is
similar to his, and indeed only claims to be somewhat shorter and more
explicit than his recent preprint~\cite{wilson:ueg}.

The reader is directed to~\cite{harpe:uniform} for a survey on uniform
growth of groups.

%%%%%%%%%%%%%%%%%%%%%%%%%%%%%%%%%%%%%%%%%%%%%%%%%%%%%%%%%%%%%%%%
\section{A Group of Non-Uniformly Exponential Growth}
\begin{figure}
  \begin{center}
    \psfrag{1}{$\9$}
    \psfrag{2}{$\8$}
    \psfrag{3}{$\mathbb 3$}
    \psfrag{4}{$\mathbb 4$}
    \psfrag{5}{$\mathbb 5$}
    \psfrag{6}{$\mathbb 6$}
    \psfrag{7}{$\mathbb 7$}
    \psfrag{x}{$x$}
    \psfrag{y}{$y$}
    \psfrag{z}{$z$}
    \hbox{\epsfig{file=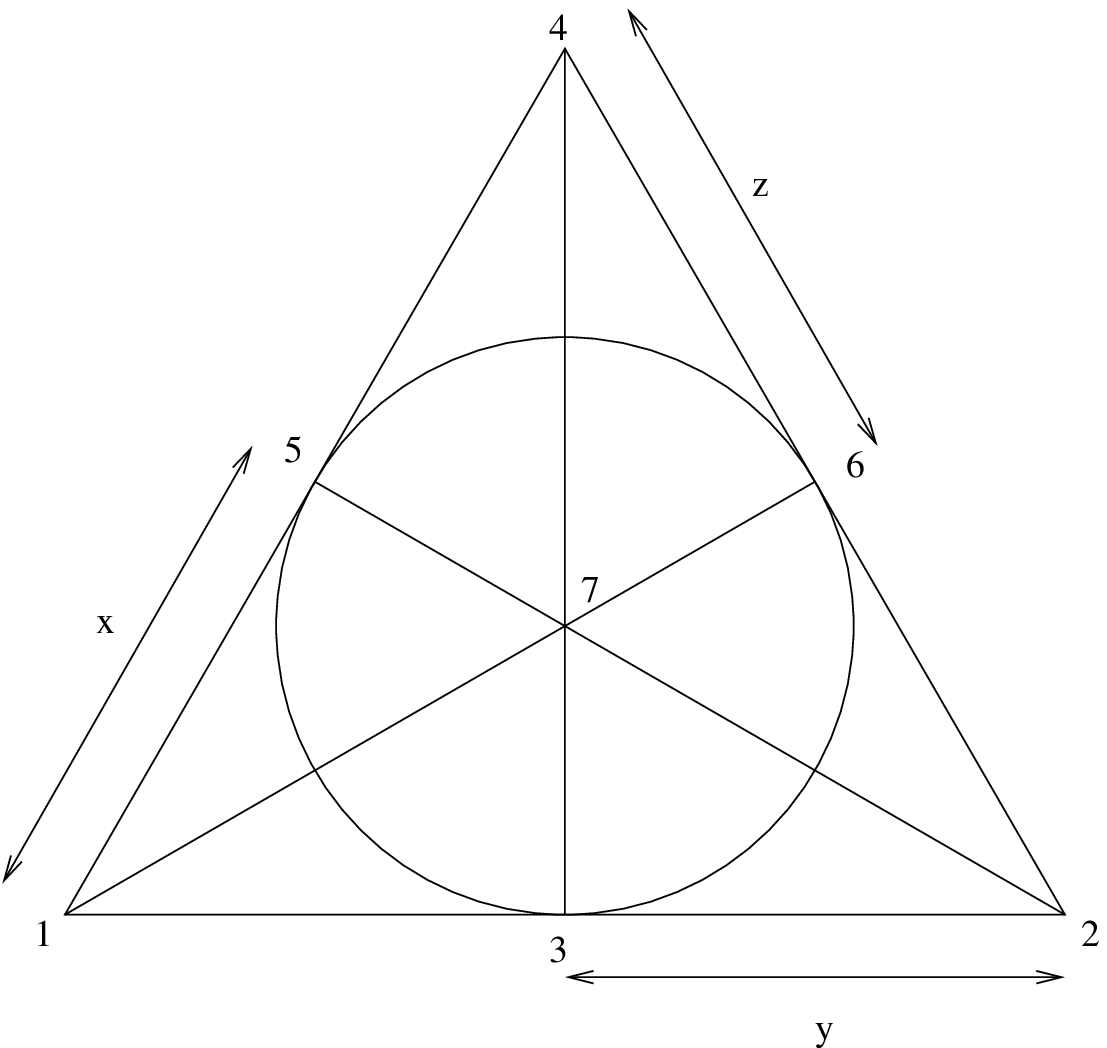,width=2.5in}\quad
      \psfrag{0}{\footnotesize $0$}
      \psfrag{0.02}{\footnotesize $0.02$}
      \psfrag{0.04}{\footnotesize $0.04$}
      \psfrag{0.06}{\footnotesize $0.06$}
      \psfrag{0.08}{\footnotesize $0.08$}
      \psfrag{0.1}{\footnotesize $0.1$}
      \psfrag{1}{\footnotesize $1$}
      \psfrag{1.2}{\footnotesize $1.2$}
      \psfrag{1.4}{\footnotesize $1.4$}
      \psfrag{1.6}{\footnotesize $1.6$}
      \psfrag{1.8}{\footnotesize $1.8$}
      \psfrag{2}{\footnotesize $2$}
      \psfrag{eta}{$\eta$}
      \begin{picture}(180,290)
        \put(0,0){\epsfig{file=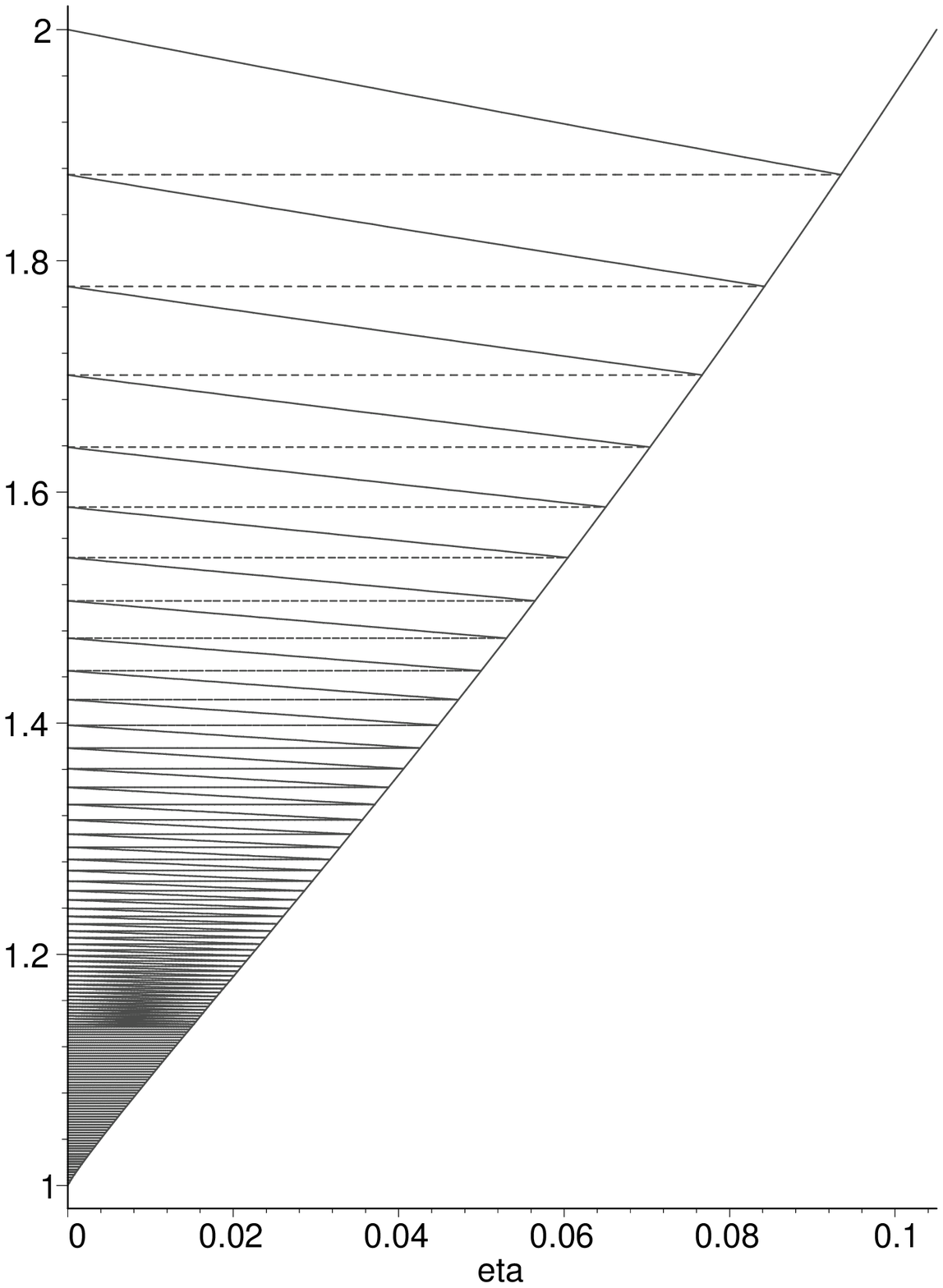,width=2.5in}}
        \put(60,240){$2^{1-\eta}$}
        \put(130,160){$\frac{30^\eta}{\eta^\eta(1-\eta)^{1-\eta}}$}
      \end{picture}}
  \end{center}
  \caption{The projective plane over $\F[2]$, and the functions in the
    proof of Proposition~\ref{prop:lowergrowth}}\label{fig:fano}
\end{figure}

First consider the group $A=\PSL(3,2)$ acting on the $7$-point
projective plane $P$ over $\F[2]$. The group $A$ is generated by $3$
reflections $x,y,z$ in $P$, as in Figure~\ref{fig:fano}.

In symbolic notation, we write
\[x=(\mathbb1,\mathbb5)(\mathbb3,\mathbb7),\quad
y=(\mathbb2,\mathbb3)(\mathbb6,\mathbb7),\quad
z=(\mathbb4,\mathbb6)(\mathbb5,\mathbb7).\]

Next we let $A$ act on the sequences $P^*$ over $P$ by
\[(p_1p_2\dots p_n)a = (p_1a)p_2\dots p_n,\]
and we let an isomorphic copy $\overline A$ of $A$ act on $P^*$ by
\[(p_1\dots p_mp_{m+1}\dots p_n)\overline a=p_1\dots
p_m(p_{m+1}a)p_{m+2}\dots p_n,
\]
where $p_1=\dots=p_{m-1}=\9$, $p_m=\8$. Alternatively, given
$g_1,\dots,g_7\in W$ and $a\in A$, we introduce the notation
$g=\pair<g_1,\dots,g_7>a$ for the permutation $g$ of $P^*$ defined by
\[(p_1p_2\dots p_n)g=(p_1a)((p_2\dots p_n)g_{p_1}).\]
In this notation, we have $\overline a=\pair<\overline
a,a,1,1,1,1,1>$.

We then define the group $W=\langle A,\overline A\rangle$ by its
action on $P^*$.

Let $G$ be a group generated by a finite generating set $S$. Its
\emph{growth rate} is
\[\lambda(G,S)=\lim_{n\to\infty}\sqrt[n]{\#B_{G,S}(n)},\]
where $B_{G,S}(n)=\setsuch{g\in G}{g=s_1\dots s_n\text{ for some
  }s_i\in S}$ is the ball of radius $n$ in $G$, with the word metric
induced by $S$. (This limit exists because $\log\#B_{G,S}(n)$ is a
subadditive function.)
  
The group $G$ has \emph{exponential growth} if $\lambda(G,S)>1$ for
one, or equivalently for any, generating set, and has
\emph{subexponential growth} otherwise. If $\lambda(G,S)=1$ and
$\#B_{G,S}(n)$ is not bounded by any polynomial function of $n$, then
$G$ has \emph{intermediate growth}. It is non-trivial to construct
groups of intermediate growth, and the first example was produced by
Grigorchuk~\cite{grigorchuk:growth} in 1983.

Note that if $\lambda_{G,S}>1$, then there exist other generating sets
$S'$ for $G$ with $\lambda(G,S')$ arbitrarily large --- for instance,
$\lambda(G,B_{G,S}(k))=\lambda(G,S)^k$. On the other hand, it is not
obvious that $\lambda(G,S)$ can be made arbitrarily close to $1$.

The group $G$ has \emph{uniformly exponential growth} if
$\inf_{\text{finite }S}\lambda(G,S)>1$.

Note that free groups, and more generally hyperbolic groups, have
uniformly exponential growth as soon as they have exponential growth.
Solvable groups~\cite{osin:entropy}, and linear
groups~\cite{eskin-m-o:uniform} in characteristic $0$, also have
uniformly exponential growth as soon as they have exponential growth.

Mikhael Gromov asked in 1981 whether there exist groups of
exponential, but non-uniformly exponential
growth~\cite[Remarque~5.12]{gromov:metriques}. This was answered
positively by John Wilson~\cite{wilson:ueg}, and is the main result of
this note:
\begin{thm}\label{thm:main}
  $W$ is a group of exponential growth, but not of uniformly
  exponential growth.
\end{thm}

The proof relies on the following propositions:
\begin{prop}\label{prop:decomp}
  $W$ satisfies the decomposition\footnote{By $W\wr A$ we mean the
    wreath product $\{f:P\to W\}\rtimes A$.} $W=W\wr A$.
\end{prop}

\begin{prop}\label{prop:free}
  $W$ contains a free monoid on $2$ generators.
\end{prop}

Given a triple $\{a,b,c\}$ of involutions acting on $P^*$, we define a
new triple $\{a',b',c'\}$ of involutions acting on $P^*$ by
\[a'=\pair<1,1,1,a,1,1,1>x,\quad b'=\pair<b,1,1,1,1,1,1>y,\quad
  c'=\pair<1,c,1,1,1,1,1>z.
\]
\begin{prop}\label{prop:extend}
  If $G$ is a perfect group generated by $3$ involutions $a,b,c$, then
  $\{a',b',c'\}$ generates $G\wr A$.
\end{prop}

\begin{prop}\label{prop:invol}
  $W$ is generated by $3$ involutions.
\end{prop}

\begin{prop}\label{prop:lowergrowth}
  Let $G$ be generated by a triple of involutions $S=\{a,b,c\}$, and
  set $S=\{a',b',c'\}$ and $H=\langle S'\rangle$. Then
  \begin{equation}\tag{$\ddagger$}\label{eq:growth}
    \lambda(H,S')\le\inf_{\eta\in(0,1)}\max\left\{\lambda(G,S)^{1-\eta},\frac{30^\eta}{\eta^\eta(1-\eta)^{1-\eta}}\right\}.
  \end{equation}
\end{prop}

\begin{proof}[Proof of Theorem~\ref{thm:main}]
  $W$ has exponential growth by Proposition~\ref{prop:free}.
  
  Pick by Proposition~\ref{prop:invol} a generating set
  $S_1=\{a,b,c\}$ of $W$ consisting of involutions; then
  $\lambda(G,S_1)\le 2$. For all $n\ge1$ apply
  Propositions~\ref{prop:lowergrowth} and~\ref{prop:extend} to
  $(W,S_n)$ to obtain $(W,S_{n+1})$. Define inductively $\Lambda_1=2$,
  and $\Lambda_{n+1}$ by solving for the unique $\eta_n\in(0,1)$ such
  that
  \begin{equation}\label{eq:main}\tag{*}
    \Lambda_{n+1}=\Lambda_n^{1-\eta_n}=30^{\eta_n}\eta_n^{-\eta_n}(1-\eta_n)^{\eta_n-1}.
  \end{equation}
  It is clear that $1<\Lambda_{n+1}<\Lambda_n$, so
  $\Lambda=\lim_{n\to\infty}\Lambda_n$ and
  $\eta=\lim_{n\to\infty}\eta_n$ exist. From~\eqref{eq:main} we have
  $\Lambda=\Lambda^{1-\eta}$, so either $\Lambda=1$ or $\eta=0$, which
  again implies $\Lambda=1$. Since $\lambda(W,S_n)\le\Lambda_n$ for
  all $n\in\N$, we have $\lim_{n\to\infty}\lambda(W,S_n)=1$.
\end{proof}

%%%%%%%%%%%%%%%%%%%%%%%%%%%%%%%%%%%%%%%%%%%%%%%%%%%%%%%%%%%%%%%%
\section{A Group Of Intermediate Growth}
Consider next the set $\tilde S=\{\tilde x,\tilde y,\tilde z\}$ of
transformations of $P^*$ defined inductively by
\[\tilde x=\pair<1,1,1,\tilde x,1,1,1>x,\quad
\tilde y=\pair<\tilde y,1,1,1,1,1,1>y,\quad \tilde z=\pair<1,\tilde
z,1,1,1,1,1>z,
\]
and consider the group $V=\langle\tilde S\rangle$.
\begin{thm}\label{thm:interm}
  $V$ is a group of intermediate growth.
  
  $V$ is locally isomorphic to $W$, in that for any $R\in\N$, there is
  $n\in\N$ such that $B_{V,\tilde S}(R)$ and $B_{W,S_n}(R)$ are
  isomorphic graphs, as seen as subsets of their respective group's
  Cayley graph.
\end{thm}
(Note that $V$ is not perfect; indeed $V/V'\cong(\Z/2)^3$. Hence $V$
does not decompose as a wreath product like $W$.)
\begin{proof}
  Apply Proposition~\ref{prop:lowergrowth} to $(V,\tilde S)$ to obtain
  $(V,\tilde S')$, and notice $\tilde S'=\tilde S$; hence
  $\lambda(V,\tilde S)=\lambda(V',\tilde S')$, so $\lambda(V,\tilde
  S)=1$ by~\eqref{eq:growth}.
  
  The groups $V$ and $W$ are contracting, i.e.\ there are constants
  $\rho=\frac12$ and $M=1$ such that for $G\in\{V,W\}$ the
  decomposition of $B(G,R)$ is a subset of $B(G,\rho R+M)\wr A$ for
  all $R\in\N$.

  Pick now $R\in\N$. There exists therefore $n\in\N$ such that the
  $n$-fold decomposition of $B(G,R)$ is a subset of $(\dots(B(G,1)\wr
  A)\wr A\dots\wr A)$. Since the generators $\tilde S$ and $S_1$ agree
  on a ball of radius $1$, this implies that the generators $\tilde S$
  and $S_n$ agree on a ball of radius $R$.
\end{proof}

%%%%%%%%%%%%%%%%%%%%%%%%%%%%%%%%%%%%%%%%%%%%%%%%%%%%%%%%%%%%%%%%
\section{Proofs}
We will use repeatedly the following facts on $A$: it has order $168$,
and is simple, hence perfect. It is generated by $\{x,y,z\}$, and also
by $\{xy,yz,zx\}$.

\begin{proof}[Proof of Proposition~\ref{prop:decomp}]
  Since $A$ acts $2$-transitively on $P$, there is $u\in A$ that fixes
  $\9$ and moves $\8$ to another point, and $v\in A$ that fixes $\8$
  and moves $\9$ to another point.
  
  Then $W$ contains $[\overline a,\overline b^u]=\pair<[\overline
  a,\overline b],1,\dots,1>$ for any $a,b\in A$, and since $A$ is
  perfect $W$ contains $\pair<\overline A,1,\dots,1>$. Similarly, $W$
  contains $[\overline a,\overline b^v]=\pair<1,[a,b],1,\dots,1>$ for
  any $a,b\in A$, so $W$ contains $\pair<1,A,1,\dots,1>$. Combining
  these, $W$ contains $\pair<W,\dots,W>$ and $A$, so $W$ contains
  $W\wr A$. The converse inclusion is obvious.
\end{proof}

\begin{proof}[Proof of Proposition~\ref{prop:free}]
  Pick $u\neq v\in A$ such that $\9u=\9v=\8$ and $\8u=\8v=\9$.
  Consider the elements $a=\overline uu,b=\overline uv,c=\overline
  vu,d=\overline vv$.  They admit the decompositions
  \begin{xalignat*}{2}
    a&=\pair<\overline u,u,1,\dots,1>(\9,\8)\sigma,&
    b&=\pair<\overline u,u,1,\dots,1>(\9,\8)\tau,\tag{$\dagger$}\label{eq:free:decomp}\\
    c&=\pair<\overline v,v,1,\dots,1>(\9,\8)\sigma,&
    d&=\pair<\overline v,v,1,\dots,1>(\9,\8)\tau,
  \end{xalignat*}
  for some permutations $\sigma\neq\tau$ of $P\setminus\{\9,\8\}$. I
  claim that $M=\{a,d\}^*$ is a free monoid; actually, we will show
  something slightly stronger, namely that $\{a,b,c,d\}^*/(a=b,c=d)$
  is freely generated by $\{a,d\}$.
  
  Consider two words $X,Y$ over $a,b,c,d$, that are not equivalent
  under $(a=b,c=d)$; we will prove by induction on $|X|+|Y|$ that they
  act differently on $P^*$. We may assume that $X$ and $Y$ are both
  non-empty, and that $X$ starts by $a$ or $b$, and $Y$ starts by $c$
  or $d$. If $|X|=|Y|=1$ then the
  decompositions~\eqref{eq:free:decomp} show that $X$ and $Y$ act
  differently on $P^*$. Otherwise, we have $|X|\equiv|Y|\mod 2$, by
  considering the action on $\9$; and furthermore we may assume
  $|X|\equiv|Y|\equiv0\mod 2$, by multiplying both $X$ and $Y$ by $a$
  on the right. Consider the decompositions
  \[X=\pair<X_1,\dots,X_7>\alpha,\quad Y=\pair<Y_1,\dots,Y_7>\beta.\]
  Then $X_1,Y_1\in\{a,b,c,d\}^*$, and $X_1$ starts by $a$ or $b$, and
  $Y_1$ starts by $c$ or $d$. We have $|X_1|=|X|/2$ and $|Y_1|=|Y|/2$,
  so by induction $X_1$ and $Y_1$ act differently on $P^*$, and hence
  so do $X$ and $Y$.
\end{proof}

\begin{proof}[Proof of Proposition~\ref{prop:extend}]
  Set $H=\langle a',b',c'\rangle$. Then $H$ contains
  $(a'b'c'b')^3=\pair<1,1,ac,ac,1,1,ca>$,
  $(b'c'a'c')^3=\pair<ba,1,1,1,1,ba,ab>$ and
  $(c'a'b'a')^3=\pair<1,cb,1,1,cb,1,bc>$; so $H$ contains
  \[v=[(a'b'c'b')^3,(b'c'a'c')^3]=\pair<1,\dots,1,[ca,ac]>\neq1.\]
  Now $\langle ac,cb,ba\rangle=G$ because $G$ is perfect, so we may
  conjugate $v$ by $(b'c'a'c')^3$ etc.\ to see that $H$ contains
  $\pair<G,\dots,G>$.
  
  Therefore $H$ contains $x=\pair<1,1,1,a,1,1,1>a'$, and similarly $y$
  and $z$, so $H=G\wr A$.
\end{proof}

\begin{proof}[Proof of Proposition~\ref{prop:invol}]
  Define $a,b,c\in W$ by
  \[a=\pair<1,\overline x,1,x,1,1,1>x,\quad b=\pair<y,1,1,\overline
  y,1,1,1>y,\quad c=\pair<\overline z,z,1,1,1,1,1>z.
  \]
  Then $(ab)^4=\pair<1,\overline x,\overline x,(x\overline
  y)^4,1,\overline x,\overline x>=\pair<1,\overline x,\overline
  x,1,1,\overline x,\overline x>$, and similarly
  $(bc)^4=\pair<1,1,1,\overline y,\overline y,\overline y,\overline
  y>$ and $(ca)^4=\pair<\overline z,1,\overline z,1,\overline
  z,1,\overline z>$.
  
  Therefore $G=\langle a,b,c\rangle$ contains
  \[u=[[(ab)^4,(bc)^4],(ca)^4]=\pair<1,\dots,1,[[\overline x,\overline
  y],\overline z]>\neq1,
  \]
  so $G$ contains all of $u$'s conjugates by $(ab)^4,(bc)^4,(ca)^4$,
  and since $A$ is perfect $G$ contains $\pair<1,\dots,1,\overline
  A>$; conjugating by $a,b,c$, we see that $G$ contains
  $\pair<\overline A,\dots,\overline A>$.
  
  Next, $G$ contains $a'=\pair<1,\overline
  x,1,1,1,1,1>a=\pair<1,1,1,x,1,1,1>x$, and similarly
  $b'=\pair<y,1,1,1,1,1,1>y$ and $c'=\pair<1,z,1,1,1,1,1>z$, so by
  Proposition~\ref{prop:extend} $G$ contains $A\wr A$; therefore
  $W=G=\langle a,b,c\rangle$.
\end{proof}

\begin{proof}[Proof of Proposition~\ref{prop:lowergrowth}]
  Consider a word $w\in\{a',b',c'\}^*$ representing an element in $H$,
  and compute its decomposition $\pair<w_1,\dots,w_7>\sigma$.  Each of
  the $w_i$'s is a word over $\{a,b,c\}$, and the total length of the
  $w_i$ is at most the length of $w$, since each $a',b',c'$
  contributes a single $a,b,c$-letter to one of the $w_i$'s.
  
  A \emph{reduced} word is a word with no two identical consecutive
  letters; we shall always assume the words we consider are reduced.
  Therefore all $aa$-, $bb$- and $cc$-subwords of the $w_i$'s should
  be cancelled; and such subwords appear in a $w_i$ whenever $w$ has a
  subword belonging to
  \[\Delta=\{a'b'a',b'c'b',c'a'c',a'c'b'a'c'a'b'c'a',b'a'c'b'a'b'c'a'b',c'b'a'c'b'c'a'b'c'\};\]
  indeed
  \begin{gather*}
    a'b'a'=\pair<1,1,1,1,b,1,1>xyx,\\
    a'c'b'a'c'a'b'c'a'=\pair<a,cb,1,1,bc,a,c>yzxzy.
  \end{gather*}
  
  \begin{lem}
    For any $n\in\N$, there are at most $30$ reduced words $w$ of
    length $n$ that contain no subword belonging to $\Delta$.
  \end{lem}
  \begin{proof}
    If $w$ contains $a'c'a'$, $b'a'b'$ or $c'b'c'$ as a subword, then
    this subword occurs either among the first $5$ or the last $5$
    letters of $w$, and $w$ is a subword of $(xyz)^\infty
    y(zyx)^\infty$, where $x,y,z$ is a cyclic permutation of
    $a',b',c'$. This gives $24$ possibilities: $3$ for the choice of
    the cyclic permutation and $8$ for the position of $zyz$ subword
    in $w$.
    
    If $w$ does not contain any such subword, then $w$ must be a
    subword of $(xyz)^\infty$ or $(zyx)^\infty$, and this gives $6$
    possibilities: $3$ for the choice of the cyclic permutation and $2$
    for the choice of $xyz$ or $zyx$.
  \end{proof}

  Fix now for every $h\in H$ a word $w_h$ of minimal length
  representing $h$; and for all $n\in\N$ let $F_n$ denote the set of
  such words of length $n$. We wish to estimate $\#F_n$.

  For any $\eta\in(0,1)$, define the following sets:
  \begin{align*}
    F_n^{>\eta}&=\setsuch{w\in F_n}{w\text{ contains at least }\eta n\text{
        subwords belonging to }\Delta},\\
    F_n^{<\eta}&=\setsuch{w\in F_n}{w\text{ contains at most }\eta n\text{
        subwords belonging to }\Delta}.
  \end{align*}
  Note that any $w\in F_n^{<\eta}$ factors as a product of at most
  $\eta n$ pieces $u_1\dots u_m$, where each $u_i$ does not contain
  any subword from $\Delta$. We therefore have
  \[\#F_n^{<\eta}\le \eta n30^{\eta n}\binom{n}{\eta n};\]
  the $\eta n$ accounting for all numbers of pieces between $1$ and
  $\eta n$, the $30^{\eta n}$ counting (according to the Lemma) the
  possible choices of each piece, and the binomial term counting for
  the respective lengths of the pieces --- they are determined by a
  selection of $\eta n$ ``separation points'' among $n$.

  Estimating the binomial coefficient $\binom{n}{\eta
    n}\approx\left(\eta^\eta(1-\eta)^{1-\eta}\right)^{-n}$, we get
  \[\lim_{n\to\infty}\sqrt[n]{\#F_n^{<\eta}}\le
  \frac{30^\eta}{\eta^\eta(1-\eta)^{1-\eta}}.\]
  
  Next, every $w\in W_n^{>\eta}$ gives by decomposition $7$ words
  $w_1,\dots,w_7$ of total length at most $(1-\eta)n$, after
  cancellation of the $aa$-, $bb$- and $cc$-subwords.
  
  For any $\epsilon>0$, there is a constant $K$ such that
  $\#B_{G,S}(n)\le K(\lambda(G,S)+\epsilon)^n$ for all
  $n\in\N$. Therefore
  \[\#F_n^{>\eta}\le\#A\binom{n+7}{7}K^7(\lambda(G,S)+\epsilon)^{(1-\eta)n};\]
  the binomial term majoring all possible partitions of the total
  length of the $w_i$'s in $7$ parts, and the other terms counting the
  number of values the $w_i$'s may assume. It follows that
  $\lim_{n\to\infty}\sqrt[n]{\#F_n^{>\eta}}\le(\lambda(G,S)+\epsilon)^{1-\eta}$
  for all $\epsilon>0$, and therefore
  \[\lim_{n\to\infty}\sqrt[n]{\#F_n^{>\eta}}\le\lambda(G,S)^{1-\eta}.\]
  
  Now $\#B_{H,S'}(n)\le n\left(\#F_n^{<\eta}+\#F_n^{>\eta}\right)$,
  and
  \begin{align*}
    \lambda(H,S')&=\lim_{n\to\infty}\sqrt[n]{\#B_{H,S'}(n)}
    \le\lim_{n\to\infty}\sqrt[n]{\#F_n^{<\eta}+\#F_n^{>\eta}}\\
    &\le\lim_{n\to\infty}\max\left\{\sqrt[n]{\#F_n^{<\eta}},\sqrt[n]{\#F_n^{>\eta}}\right\}\\
    &\le\max\left\{\lambda(G,S)^{1-\eta},\frac{30^\eta}{\eta^\eta(1-\eta)^{1-\eta}}\right\}.
  \end{align*}
\end{proof}

%%%%%%%%%%%%%%%%%%%%%%%%%%%%%%%%%%%%%%%%%%%%%%%%%%%%%%%%%%%%%%%%
\bibliography{mrabbrev,people,math,grigorchuk,bartholdi}
\end{document}